\documentclass[12pt]{article}
\usepackage[final]{epsfig}
\usepackage{amsfonts, amsmath, amssymb, amsthm}
\usepackage{latexsym}
\usepackage{color}
\usepackage[T1]{fontenc}

\numberwithin{equation}{section}

\newtheorem{theorem}{Theorem}[section]
\newtheorem{lemma}{Lemma}[section]
\newtheorem{example}{Example}[section]
\newtheorem{proposition}{Proposition}[section]
\newtheorem{definition}{Definition}[section]
\newtheorem{corollary}{Corollary}[section]

\newcommand{\R}{\mathbb{R}} 

\newcommand{\Z}{\mathbb{Z}}
\newcommand{\N}{\mathbb{N}}

\newcommand{\be}{\begin{equation}}
\newcommand{\ee}{\end{equation}}
\newcommand{\calP}{{\mathcal{P}}}
\newcommand{\calPoneN}{{\mathcal{P}}_{\frac{1}{N}}}
\newcommand{\calM}{{\mathcal{M}}}
\newcommand{\calI}{{\mathcal{I}}}
\newcommand{\eps}{{\varepsilon}}

\newcommand{\lambdabar}{{\overline{\lambda}}}

\begin{document}

\begin{center}
\begin{LARGE}
   Breaking the curse of dimension \\
  in multi-marginal Kantorovich optimal transport \\
       on finite state spaces \\
\end{LARGE}

\normalsize
\vspace{0.4in}

Gero Friesecke$^1$ and Daniela V\"ogler$^1$ \\[1mm]
$^1\,$Faculty of Mathematics, Technische Universit\"at M\"unchen  \\
gf@ma.tum.de, voegler@ma.tum.de \\[2mm]
 
\end{center}


{\bf Abstract.} We present a new ansatz space for the general symmetric multi-marginal Kantorovich optimal transport problem on finite state spaces which reduces the number of unknowns from 
$\tbinom{N+\ell-1}{\ell-1}$ to $\mbox{\small $\ell\cdot(N+1) $}$, where $\ell$ is the number of marginal states and $N$ the number of marginals. 

The new ansatz space is a careful low-dimensional enlargement of the Monge class, which corresponds to $\mbox{\small $\ell\cdot(N-1) $}$ unknowns, and cures the insufficiency of the Monge ansatz, i.e. we show that the Kantorovich problem always admits a minimizer in the enlarged class, for arbitrary cost functions. 

Our results apply, in particular, to the discretization of multi-marginal optimal transport with Coulomb cost in three dimensions, which has received much recent interest due to its emergence as the strongly correlated limit of Hohenberg-Kohn density functional theory. In this context $N$ corresponds to the number of particles, motivating the interest in large $N$. 

\section{Introduction} \label{sec:Intro} 
In this paper we study the multi-marginal optimal transport problem 
\be \label{kant}
   \mbox{Minimize }C[\gamma] = \int_{X^N} c_N(x_1,...,x_N) \, d\gamma(x_1,...,x_N) \mbox{ over }\gamma\in\calP_{sym}(X^{N}) 
   \mbox{ subject to }\gamma\mapsto\lambda_* 
\ee
for finite state spaces  
\be \label{X}
    X = \{a_1,...,a_\ell \}
\ee
consisting of $\ell$ distinct points. 
Here $\calP_{sym}(X^{N})$ denotes the set of symmetric probability measures on $X^N$, where symmetric means
\be \label{sym}
   \gamma(A_1\times\cdots\times A_N) = \gamma(A_{\sigma(1)}\times\cdots\times A_{\sigma(N)}) 
   \mbox{ for all subsets } A_1,...,A_N \mbox{ of } X 
   \mbox{ and all permutations } \sigma, 
\ee
$c_N \, : \, X^N\to\R\cup\{+\infty\}$ is an arbitrary cost function, 
$\lambda_*$ denotes a given probability measure on $X$, and the notation $\gamma\mapsto\lambda_*$ means that $\gamma$ has equal one-point-marginals $\lambda_*$, i.e.
\be \label{marginal}
   \gamma(X^{i-1}\times A_i\times X^{N-i}) = \lambda_*(A_i) \mbox{ for all subsets }A_i\mbox{ of }X 
   \mbox{ and all }i=1,...,N. 
\ee
   
Multi-marginal problems \eqref{kant} arise in economics \cite{CE10, CMN10}, electronic structure \cite{CFK11, BDG12}, image processing \cite{AC11, RPDB12}, mathematical finance \cite{BHP13, GHT14}, and optimal assignment problems \cite{Pi68, GH88, Po94, APS99, Pe11}; in some cases finite state spaces \eqref{X} appear directly, and in others they play the role of natural discretizations of continuous spaces (see e.g. \cite{AC11, CFM14}). 

A prototypical example of \eqref{kant}, \eqref{X} is discretized multi-marginal optimal transport with Coulomb cost in $\R^3$, where $X$ is a collection of $\ell$ discretization points 
$a_i$ in $\R^3$ and
\be \label{coul}
    c_N(x_1,...,x_N) = \sum_{1\le i < j \le N} \frac{1}{|x_i-x_j|},
\ee
with $|\, \cdot \, |$ being the euclidean distance. 
%
Physically, the $a_i$ correspond to possible sites and the Kantorovich problem \eqref{kant} corresponds to optimally assigning the positions of $N$ particles to the available sites subject to the constraint that all sites must be occupied according to the prescribed marginal measure. (Even in the prototype case of uniform marginal, i.e. when each site must be occupied equally often, such a task may require stochastic superposition of configurations. Readers who have the contrary impression are advised to have a look at Example \ref{E:counter}.) 
In its continuous form with $X=\R^3$, the problem \eqref{kant}, \eqref{coul} 
has received much interest recently in both mathematics and physics due to its emergence as a strongly correlated limit of Hohenberg-Kohn density functional theory. In this context the restriction to symmetric plans in \eqref{kant} is inherited from the symmetries of the quantum problem, and the optimal cost as a functional of the marginal measure -- or single-particle density -- is known as the SCE (strictly correlated electrons) functional. See \cite{Se99, SGS07} for the original formulation and derivation of the limit problem in the physics literature, \cite{CFK11, BDG12} for its reformulation as an OT problem within a rigorous functional-analytic setting, and \cite{CFK11, BD17, CFK17} for the rigorous justification as $\Gamma$-limit of the Hohenberg-Kohn functional for, respectively, $N=2$, $N=3$, and general $N$. An analogous general-$N$ result for a ``relaxation'' of the HK functional can be found in \cite{Le17}. 

The prototypical marginal on finite state spaces \eqref{X} is the uniform marginal, i.e. 
$\lambda_*=\lambdabar$ where $\lambdabar$ denotes the uniform probability measure on $X$,
\be \label{unif}
    \lambdabar = \sum_{i=1}^\ell \tfrac{1}{\ell} \, \delta_{a_i}.
\ee
Here and below, $\delta_{a_i}$ denotes the Dirac measure on the point $a_i$. This marginal measure arises directly in assignment problems, and via equi-mass discretization \cite{CFM14} from continuous problems: given any absolutely continuous probability measure on $\R^d$, divide $\R^d$ into regions carrying equal mass and chose the $a_i$ to be representative points for each region. 
The uniform marginal is the natural discrete analogue of the absolutely continuous marginals on euclidean space, because it admits Monge states
\be \label{I:Monge}
   \gamma = \sum_{\nu=1}^\ell \tfrac{1}{\ell} 
            \delta_{T_1(a_\nu)} \otimes \cdots \otimes \delta_{T_N(a_\nu)} 
            \mbox{ for $N$ permutations } T_1,...,T_N\, : \, X\to X.
\ee 
Note that the requirement that a map $T \, : \, X\to X$ be a permutation, i.e. that 
$T(a_\nu)=a_{\tau(\nu)}$ for all $\nu$ and some permutation $\tau \, : \, \{1,...,\ell\}\to \{1,...,\ell\}$ of indices, 
is the same as requiring that the map pushes the uniform measure forward to itself. In the present setting of symmetric multi-marginal problems, the ansatz must be trivially adapted to the symmetrization 
\be \label{I:Monge'}
   \gamma' = S \, \gamma, \; \gamma \mbox{ as in \eqref{I:Monge}}, 
\ee
where $S$ is the symmetrization operator in $N$ variables (see \eqref{S}). Restricting the minimization in \eqref{kant} to states of form \eqref{I:Monge'} gives the corresponding Monge OT problem. (When the cost $c_N$ is symmetric, i.e. invariant under permuting its arguments, as is the case in many examples of interest such as \eqref{coul}, the requirement that $\gamma$ be symmetric can be dropped in the Kantorovich problem \eqref{kant} -- and likewise the symmetrization operator $S$ can be  omitted in the Monge problem -- without altering the minimum cost.) 

In the present finite-dimensional context it is clear that the Monge ansatz entails a spectacular reduction of computational complexity. It reduces the number of unknowns from combinatorial in both $N$ and $\ell$ (see Theorem \ref{T:extreme'} for the precise numbers) to only $\ell\cdot (N-1)$,  because each map $T_k$ is specified by its $\ell$ values $T_k(a_1),...,T_k(a_\ell)$ and one may assume $T_1=id$, by re-ordering the sum in \eqref{I:Monge}. Thus unlike the Kantorovich problem, the Monge problem remains computationally feasible for large $\ell$ and $N$.

Unfortunately the Monge ansatz is not always sufficient to obtain the optimal cost when $N\ge 3$.
(For $N=2$, i.e. two marginals, it suffices thanks to the celebrated Birkhoff-von Neumann theorem \cite{Bi46, vN53}; optimal transportation theory (see e.g. \cite{Vi09}) provides the analogous result for continuous two-marginal problems under very general conditions.) For continuous problems with $N\ge 3$, whether the ansatz \eqref{I:Monge} works appears to depend on subtle properties of the cost function and even on the ambient space dimension; for the Coulomb cost with $X=\R^3$ it is presently unknown. See \cite{GS98, He02, Ca03, Pa11, CDD13} for interesting examples where minimizers of \eqref{kant} are of Monge form, with the first result of this type appearing in a fundamental paper by Gangbo and \'{S}wi\k{e}ch \cite{GS98}. Examples of non-Monge minimizers can be found in \cite{Pa10, FMPCK13, Pa13}, and see \cite{MP17} for an example with unique non-Monge minimizer. For discrete assignment problems with finite state space $X$ and $N=3$, it is known that there exist ``non-integer vertices'' of the -- suitably renormalized -- polytope of probability measures on $X^3$ with uniform marginals \cite{Kr07, LL14}. This can be shown after some work \cite{Fi14} to imply the existence of cost functions with unique non-Monge minimizer. 
The following simple example is taken from \cite{Fr18}. 
\begin{example} \label{E:counter} The unique minimizer of the Kantorovich problem with $X$ given by three equi-spaced points (physically: sites) on the real line, i.e. $X=\{1,2,3\}\subset\R$, $N=3$ (physically: three particles), and the cost
\be \label{counterex}
   c_N(x_1,...,x_N) = \sum_{1\le i < j\le N} c(|x_i-x_j|), \;\; c(r)=(r - \tfrac{3}{4})^2
\ee
(physically: the particles are mutually connected by springs of equilibrium length $\tfrac{3}{4}$) is uniquely minimized by $\gamma_* = S \gamma$ where $\gamma=\tfrac{1}{2}(\delta_1\otimes\delta_1\otimes\delta_2 + \delta_2\otimes\delta_3\otimes\delta_3)$ and $S$ is the symmetrization operator \eqref{S}. This $\gamma_*$ is not a symmetrized Monge state. 

\end{example}
The main result of this paper (Theorem \ref{T:OTnew} a)) is that a careful low-dimensional enlargement of the Monge class, where each state requires $\ell\cdot(N+1)$ instead of $\ell\cdot(N-1)$ parameters, cures the insufficiency of the Monge ansatz, i.e. the Kantorovich problem \eqref{kant} always admits a minimizer in the enlarged class. We propose to call the states in this class {\it sparse averages of extremal states}, or SAE states for short, because they are constructed by averaging a small number of extreme points of the convex set $\calP_{sym}(X^N)$ (see Section \ref{sec:sparse}). Moreover, as in the Monge case, for pairwise costs (such as \eqref{coul}) the optimal cost in \eqref{kant} agrees with that of an explicit reduced problem which involves only two-point probability measures; that is to say all high-dimensional objects (like $\calP_{sym}(X^N)$) and operations (like integration over $X^N$) can be eliminated. See Theorem \ref{T:OTnew} b). 

We now describe the new, sufficient ansatz in the prototypical case of uniform marginal $\lambdabar$, where it can be compared to the Monge ansatz \eqref{I:Monge'}. 
\begin{enumerate}
\item Take $N$ maps $T_1,...,T_N\, : \, X\to X$. 
\item Drop the rigid weights $\tfrac{1}{\ell}$ in the Monge ansatz \eqref{I:Monge} and replace them by flexible site weights $\alpha^{(\nu)}\ge 0$ ($\nu=1,...,\ell$) which sum to $1$, and form the otherwise analogous $N$-point measure 
\be \label{I:SAE} 
    \gamma = \sum_{\nu=1}^\ell \alpha^{(\nu)} 
    S \, \delta_{T_1(a_\nu)}\otimes\cdots\otimes\delta_{T_N(a_\nu)}. 
\ee 
In particular, the weights of some sites may be zero, i.e. the $N$-point state might only use the 
values of the $T_k$ at fewer than $\ell$ sites. 
\item Drop the rigid constraints in \eqref{I:Monge} that all the $T_k$ individually preserve the uniform measure, i.e. attain each possible value exactly once. Instead, consider the empirical value distribution of the {\it ensemble of maps} $\{T_k\}_{k=1,...,N}$ at the site $a_\nu$, 
\be \label{I:valdist}
    \lambda^{(\nu)} := \frac{1}{N} \sum_{k=1}^N \delta_{T_k(a_\nu)},
\ee
and require that the {\it average over the sites $a_\nu$} of this value distribution associated with the weights $\alpha^{(\nu)}$ in \eqref{I:SAE} equals the uniform measure, i.e. impose the constraint
\be \label{I:constraint}
   \sum_{\nu=1}^\ell \alpha^{(\nu)}\lambda^{(\nu)} = \lambdabar.
\ee
\end{enumerate}
It is easy to check that the constraint \eqref{I:constraint} precisely guarantees the marginal condition $\gamma\mapsto\lambdabar$. We also emphasize that \eqref{I:SAE} is the symmetrization of a measure concentrated on a single graph over the marginal space $X$, not on several graphs as considered e.g. in \cite{MP17}. 

To summarize: 
\begin{definition} \label{D:SAE} A probability measure on $X^N$ is called an SAE state with uniform marginal if and only if it is of form \eqref{I:SAE} for some maps $T_1,...,T_N\, : \, X\to X$ such that the associated empirical value distributions \eqref{I:valdist} satisfy \eqref{I:constraint}. An SAE state with general marginal $\lambda_*$ is defined in the same manner, except that the right hand side of eq. \eqref{I:constraint} needs to be replaced by $\lambda_*$. 
\end{definition}

\begin{figure}
\includegraphics[width=0.49\textwidth]{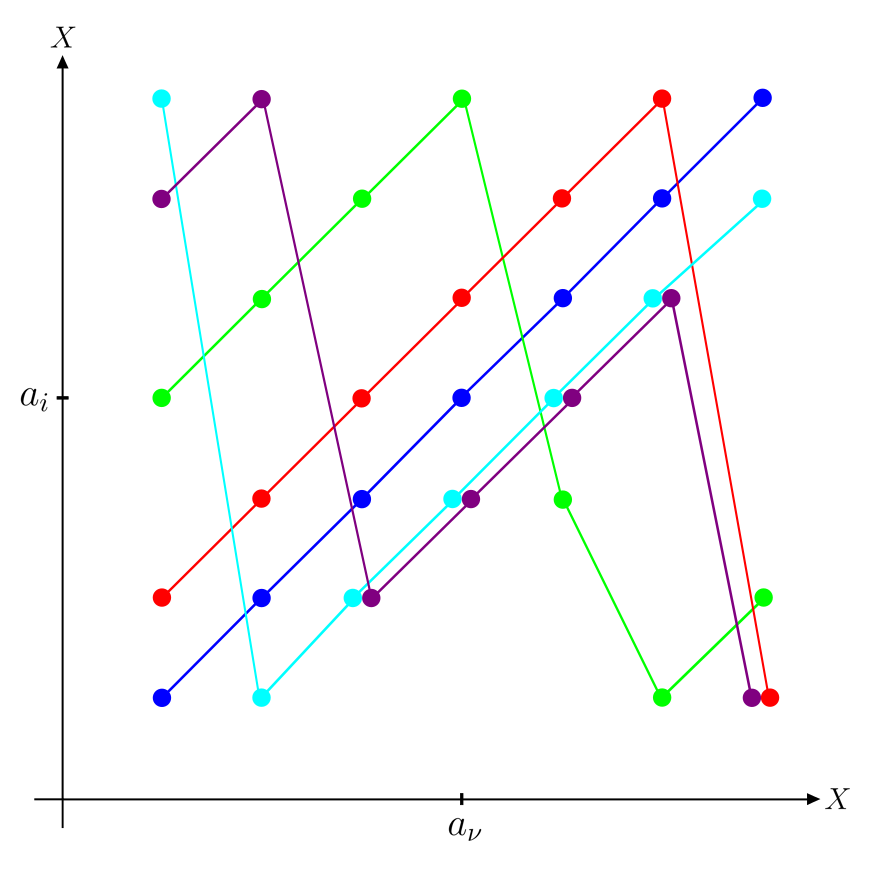}
\includegraphics[width=0.49\textwidth]{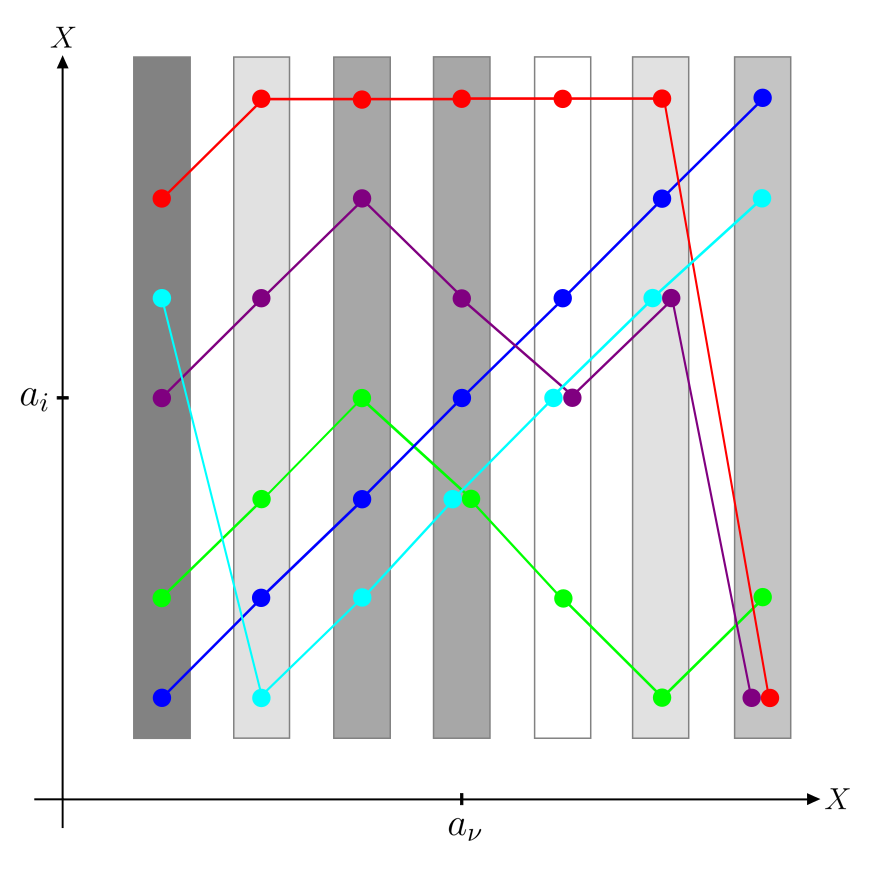}
\caption{Left: Monge state. Right: SAE state. The maps $T_k$ are indicated by different colors, and the site weights $\alpha^{(\nu)}$ by the gray level of the boxes. For the Monge state, the total value distribution (left hand side of \eqref{I:constraint}) is uniform, because each map contributes exactly one point to each site $a_i$. For the SAE state, it is also uniform even though the maps contribute unequally to the different sites. The reader may check that if the site weights (normalized to sum to 7 to correspond to point counting) are chosen as 2, 0.5, 1.5, 1.5, 0, 0.5, 1, the contribution from all maps to each site $a_i$ is still $5$ points.} \label{F:SAE}
\end{figure}

For a graphical representation of a typical SAE state see Figure \ref{F:SAE}. 

It is instructive to directly understand why each Monge state \eqref{I:Monge'} is an SAE state without passing through $N$-point measures and the computation of marginals. I.e., we would like to understand why the measure $\lambda^{av} := \sum_{\nu=1}^\ell \alpha^{(\nu)}\lambda^{(\nu)}$ is equal to $\lambdabar$ if the $\alpha^{(\nu)}$ are $\tfrac{1}{\ell}$ and the $T_k$ individually preserve the uniform measure. To see this, first change the order of averaging over maps and sites, then use constancy of the weights and the fact that the $T_k$ are measure preserving, i.e. -- in informal notation --
\begin{eqnarray*}
   \lambda_i^{av} &=& \, \mbox{average over sites}\Bigl(\mbox{average over maps}\bigl(\mbox{probability of value }a_i\bigr)\;\Bigr) \\
                  &=& \mbox{average over maps}\Bigl(\underbrace{\mbox{average over sites}\bigl(\mbox{probability of value }a_i\bigr)}_{=\, \tfrac{1}{\ell} \mbox{ since each value $a_i$ appears exactly once}.} \Bigr)
\end{eqnarray*}
This calculation reveals the extra freedom in symmetric problems not used by the Monge ansatz that the value distribution of the maps only needs to be right ``on average''. A deeper fact is that using this freedom requires flexible site weights; i.e., SAE states with equal weights $\alpha^{(\nu)}=\tfrac{1}{\ell}$ which satisfy the average condition \eqref{I:constraint}  are already Monge states. See Section \ref{sec:Monge}. 
    
Why the new ansatz, unlike \eqref{I:Monge'}, is sufficient to solve the Kantorovich problem is of course less elementary. A key step is to achieve a good understanding of the convex geometry 
of the set of $N$-representable probability measures on $X^2$ and that of symmetric probability measures on $X^N$, by using methods from convex analysis, differential geometry, and discrete mathematics. In particular we classify the extreme points of both sets, and show that the two-point-marginal map is a bijection between them. One analytically and probabilistically interesting consequence of our analysis, described in Sections \ref{sec:extremal} and \ref{sec:extremal'}, is that 
the extreme points are unique {\it maximizers} of the Wasserstein cost, respectively the Gangbo-\'{S}wi\k{e}ch cost \eqref{GS}, with respect to the discrete metric. The latter cost is the natural multi-marginal analogue of the Wasserstein cost, and was first introduced -- in the context of optimal transport on euclidean spaces -- in \cite{GS98}.

Finally we remark that the SAE ansatz space might turn out to be useful even for problems which admit Monge minimizers. It might serve either as a stepping stone to proving existence of the latter or -- because of its less rigid nature and its robustness under changes of state space and cost -- as a suitable space within which to perform iterative updates in a numerical algorithm. 
%
%

\section{Extremal $N$-representable two-point measures}
\label{sec:extremal}
Throughout this paper, $X$ denotes the finite state space \eqref{X}, and $\calP(X)$ denotes the set of probability measures on $X$. Because of the finiteness of $X$ we can identify any probability measure $\lambda\in\calP(X)$ with a vector 
$(\lambda_1,...,\lambda_\ell)\in\R^\ell$ satisfying $\lambda_j\ge 0$ and $\sum_j\lambda_j=1$, whose components are the probabilities
\be \label{lambdaj}
    \lambda_j = \lambda ( \{ a_j \} ).
\ee
We will frequently use the shorthand notation
\be \label{shorthand}
    \delta_i := \delta_{a_i}, \;\;\;
    \delta_{i_1...i_N} := \delta_{a_1}\otimes\cdots\otimes\delta_{a_N}
\ee
for the Dirac measures centered on the point $a_i\in X$ respectively $(a_{i_1},...,a_{i_N})\in X^N$. 
The following concept, introduced recently in \cite{FMPCK13}, plays an important role throughout this paper. The quantum analogue, where one is dealing with ``$N$-body density matrices'' (operators on the $N$-fold antisymmetric tensor product of a Hilbert space $X$) instead of $N$-point probability measures (measures on the $N$-fold cartesian product of a set $X$) has been introduced a long time ago (see \cite{CY02} for a textbook account). 
\begin{definition}[$N$-representability]  \label{D:Nrep}  Let $N\ge 2$. For any $k\in\{2,...,N\}$, a probability measure $\mu$ on $X^k$ is called \emph{$N$-representable} if there exists a symmetric probability measure $\gamma$ on $X^N$ such that its $k$-point marginal (see \eqref{Mk}) equals $\mu$. Any such $\gamma$ is called a \emph{representing measure} of $\mu$.  
\end{definition}
For elementary examples of measures which are not $N$-representable we refer the reader to
\cite{FMPCK13}. In the following, we denote the set of $N$-representable probability measures on $X^k$ by $\calP_{N-rep}(X^k)$. Our principal interest is in the case $k=2$, i.e. in the set 
$\calP_{N-rep}(X^2)$ of $N$-representable two-point probability measures. 

This reduced space of measures is ``dual'' to the space of pairwise cost functions
\be \label{pairwise}
  c_N(x_1,...,x_N)=\sum_{1\le i<j\le N} c(x_i,x_j) \; \mbox{ for some }
  c\, : \, X\times X\to\R\cup\{+\infty\}
\ee 
which are typical in applications. As pointed out in \cite{FMPCK13}, because of the elementary identity
\be \label{elid}
  \int_{X^N} \sum_{1\le i<j\le N} c(x_i,x_j) \, d\gamma(x_1,...,x_N) = 
  \binom{N}{2}\int_{X^2}c(x,y) \, d(M_2\gamma)(x,y)
\ee
which shows that the Kantorovich functional \eqref{kant} with pairwise cost \eqref{pairwise} only depends on the two-point marginal, we have that
\be \label{redkant}
   \min_{{\gamma\in\calP_{sym}(X^N)}\atop{\gamma\mapsto\lambda_*}} C[\gamma] = 
   \min_{{\mu\in\calP_{N-rep}(X^2)}\atop{\mu\mapsto\lambda_*}} \binom{N}{2} 
   \int_{X^2} c(x,y)\, d\mu(x,y).
\ee
That is to say the minimal cost in the high-dimensional problem on the left is the same as that of the dimension-reduced problem on the right. Eq.~\eqref{redkant} is the classical analogon of expressing the ground state energy of a quantum system with $N$ electrons via a minimization problem for $N$-representable two-body density matrices (see \cite{CY02}). 

The catch -- familiar from the quantum case -- is that the set of admissible trial states on the right is unknown, and was only constructed by using states in the high-dimensional space on the left. The only case in which this set can be understood in a straightforward manner is the ``non-multi-marginal'' case $N=2$: a probability measure $\mu$ on $X^2$ is two-representable if and only if it is symmetric. Our interest, however, motivated by the role of $N$ as the number of particles in physical systems, is in arbitrary $N$.
\\[2mm]
Thanks to the finiteness of $X$, $\calP_{N-rep}(X^2)$ is a finite-dimensional compact convex set, and therefore - by Minkowski's theorem (see e.g. \cite{Ho94}) - the convex hull of its extreme points. (Recall the following standard notions from convex analysis \cite{Ro97, Ho94}: the convex hull $conv(K)$ of a set $K$ is the set of finite convex combinations 
$$
    x = \alpha_1 x_1 + \cdots + \alpha_M x_M \mbox{ for some }M\in\N, 
    \mbox{ some }x_i\in K, \mbox{ and some } \alpha_i\ge 0 \mbox{ with }\sum_i\alpha_i = 1;
$$
and a point $x$ in a convex set $K$ is called an extreme point if, whenever $x=\alpha x_1 + (1-\alpha)x_2$ for some $x_1$, $x_2\in K$ and some $\alpha\in(0,1)$, we have that $x_1=x_2=x$.) 
\\[2mm]
The extreme points of $\calP_{N-rep}(X^2)$ can be determined explicitly. In the result below, an important role is played by a certain subset of the one-point probability measures, the $\frac{1}{N}$-quantized one-point probability measures:
\be \label{P1overN}
   \calPoneN(X) = 
   \left\{\lambda \in \calP(X) \, : \, \lambda_i\in\frac{1}{N}\Z \mbox{ for all }i\in \{1,...,\ell\} \right\}. 
\ee
\begin{theorem}[Extreme $N$-representable measures] \label{T:extreme}
A probability measure $\mu$ on $X^2$ is an extreme point of the set $\calP_{N-rep}(X^2)$ of $N$-representable two-point probability measures if and only if it is of the form 
\be \label{extr}
    \lambda \otimes \lambda + \frac{1}{N-1} \Bigl( \lambda \otimes \lambda - 
    \sum_{i=1}^\ell \lambda_i \delta_i\otimes \delta_i \Bigr)
\ee
for some $\frac{1}{N}$-quantized one-point probability measure $\lambda$ (i.e. some $\lambda\in \calPoneN(X)$).
\end{theorem}
Note that the measure \eqref{extr} has marginal $\lambda$, and contains correlations (second term in \eqref{extr}); moreover these correlations have a universal structure which depends on $N$ but not on $\lambda$. The latter fact can be expressed more concisely by introducing the following
universal marginal-to-correlated-state map $\varphi_N\, : \, \calP(X)\to \calP(X^2)$
\be \label{UN}
  \varphi_N(\lambda) := \mbox{expression } \eqref{extr}.
\ee 
Note that this map is well-defined on all of $\calP(X)$. Moreover the image measure $\varphi_N(\lambda)$ has marginal $\lambda$; in particular $\varphi_N$ is injective. Theorem \ref{T:extreme} says that the set of extreme points of $\calP_{N-rep}(X^2)$ is equal to $\varphi_N(\calPoneN(X))$, i.e. it is the image of the set of $\frac{1}{N}$-quantized one-point measures under the marginal-to-correlated state map. 
\\[2mm]
Theorem \ref{T:extreme} allows, in particular, to determine the number of extreme points.
\begin{corollary} \label{C:counting} The set $\calP_{N-rep}(X^2)$ has precisely
$\tbinom{N+\ell-1}{\ell - 1}$ extreme points. 
\end{corollary}
{\bf Proof of the corollary.} Since the set of extreme points is the image of $\calPoneN(X)$ under the injective map \eqref{UN}, it suffices to determine the cardinality of $\calPoneN(X)$. Introduce the following symbolic representation of $\lambda\in\calPoneN(X)$ by a bitstring:
$$
  x = \underbrace{0 \cdots 0}_{N\lambda_1} 1 \underbrace{0 \cdots 0}_{N\lambda_2} 1 \cdots
      \underbrace{0 \cdots 0}_{N\lambda_\ell}.
$$
Note that this bitstring consists, in total, of $(\ell-1)$ ones and $N$ zeroes, and that $\lambda$ can be recovered from $x$ via
$$
  \lambda_i = \mbox{no. of zeroes between the $(i-1)^{th}$ and $i^{th}$ one in $x$}.
$$ 
Hence the cardinality of $\calPoneN(X)$ equals that of the set of corresponding bitstrings. 
Since the latter consists of all sequences of zeroes and ones of length $N+\ell-1$ containing exactly $\ell-1$ ones, it has cardinality 
$$
     \tbinom{N+\ell-1}{\ell - 1}.
$$
This establishes the corollary. 
\\[2mm]
{\bf Proof of Theorem \ref{T:extreme}.} First we show that the set of $N$-representable two-point probability measures is equal to the convex hull of the points \eqref{extr}, i.e. that
\be \label{subset}
    \calP_{N-rep}(X^2) = conv\bigl(\varphi_N(\calPoneN(X))\bigr),
\ee
where $\varphi_N$ is the map \eqref{UN}. This statement is proved for $N=2$ in \cite{FMPCK13} and for general $N$ in \cite{Fi14}. For completeness we include a short proof of our own. We then show that each element of $\varphi_N(\calPoneN(X))$ is extreme, i.e. that for any $\lambda\in\calPoneN(X)$
\be \label{supset}
    \varphi_N(\lambda) \not\in 
    conv\Bigl(\varphi_N(\calPoneN(X)) \backslash \{\varphi_N(\lambda)\}\Bigr). 
\ee
{\bf Proof of \eqref{subset}.} This implication is more or less straightforward. For completeness we include the details. We start from the fact that, directly from the definition of $\calP(X^N)$, 
\be \label{PXN}
     \calP(X^N) = conv \bigl\{ \delta_{i_1...i_N} \, : \, i_1,...,i_N\in \{1,...,\ell\} \, \bigr\}.
\ee
Here and below we use the notation \eqref{shorthand}. 
The subset $\calP_{sym}(X^N)$ of symmetric probability measures on $X^N$ (see \eqref{sym}) 
is the image of $\calP(X^N)$ under the linear symmetrization operator $S \, : \, \calP(X^N) \to \calP(X^N)$ defined by 
\be \label{S}
   (S\gamma )(A_1 \times\cdots\times A _N) = \frac{1}{N!} \sum_{\sigma\in S_N} 
   \gamma\bigl( A_{\sigma(1)} \times \cdots \times A_{\sigma(N)}\bigr) 
   \mbox{ for all }A_1,...,A_N\subseteq X.
\ee
Here $S_N$ denotes the group of permutations $\sigma \, : \{ 1,...,N\} \to \{ 1,...,N\}$. Consequently 
\be \label{MNS}
   \calP_{sym}(X^N) = conv\{ S \delta_{i_1...i_N} \, : \,
    1 \le i_1 \le ... \le i_N\le \ell\}. 
\ee
The set $\calP_{N-rep}(X^2)$ is, by the definition of $N$-representability, the image of $\calP_{sym}(X^N)$ under the linear map $M_2 \, : \, \calP(X^N) \to \calP(X^2)$ from an $N$-point measure to its two-point marginal (see \eqref{Mk} below). It follows that 
\be \label{M2S}
    \calP_{N-rep}(X^2) = conv\{ M_2 S \delta_{i_1...i_N} \, : \,
    1 \le i_1 \le ... \le i_N\le \ell\}. 
\ee
Here and below, the map from $N$-point measures to their $k$-point marginals is denoted by $M_k$, that is to say for any $k\in\{ 1,...,N\}$ and any $\gamma\in \calP(X^N)$ we define
\be \label{Mk} 
   (M_k\gamma)(A) = \gamma(A\times X^{N-k}) \mbox{ for all }A\subseteq X^k, 
\ee
with the convention that $M_N$ is the identity. 
The $N$-representable measures appearing on the right hand side of \eqref{M2S} can be evaluated explicitly, by partitioning the sum over all permutations $\sigma$ according to the value of $\sigma(1)$ and $\sigma(2)$: 
\begin{eqnarray}
   N! \, M_2 S \delta_{i_1...i_N} & = & M_2 \sum_{\sigma\in S_N} \delta_{i_{\sigma(1)} ... i_{\sigma(N)}} \nonumber \\
   & = & M_2 \sum_{{m,n=1}\atop{m\neq n}}^N \sum_{{\sigma\in S_N}\atop{\sigma(1)=m,\, \sigma(2)=n}} \delta_{i_m i_n i_{\sigma(3)} ... i_{\sigma(N)}} \nonumber \\
   & = & (N-2)! \sum_{{m,n=1}\atop{m\neq n}}^N \delta_{i_m i_n}, \label{M2'}
\end{eqnarray}
with the factor $(N-2)!$ in the last expression arising as the number of permutations with $\sigma(1)=m$, $\sigma(2)=n$. The corresponding one-point marginal is
\be \label{M1'}
   M_1 S \delta_{i_1...i_N} = \frac{1}{N(N-1)} M_1 \sum_{{m,n=1}\atop{m\neq n}}^N \delta_{i_m i_n} = \frac{1}{N} \sum_{m=1}^N \delta_{i_m}.   
\ee
Given the explicit expressions \eqref{M2'}, \eqref{M1'}, it is now a straightforward matter to infer the following lemma.
\begin{lemma}[One- and two-point marginals of symmetrized Dirac measures] \label{L:extreme} 
We have \\
a) $\{ M_1 S\delta_{i_1...i_N} \, : \, 1 \le i_1\le ... \le i_N\le \ell\} = \calPoneN(X)$ \\
b) For any $i_1,...,i_N$ as above, $M_2 S \delta_{i_1...i_N} = \varphi_N(M_1 S \delta_{i_1...i_N})$.
\end{lemma}
Clearly, this lemma combined with \eqref{M2S} establishes \eqref{subset}, completing the proof of the ``easy part'' of Theorem \ref{T:extreme}. It remains to verify the lemma.
\\[2mm]
{\bf Proof of Lemma \ref{L:extreme}.} Given any $i_1,...,i_N\in\{ 1,...,\ell\}$, we introduce the numbers  
\be \label{lambdai}
   \lambda_i := \frac{1}{N}\sharp \bigl\{ k\in\{1,...,N\} \, : \, i_k=i\bigr\} 
   \in \left\{ 0,\tfrac{1}{N},...,\tfrac{N-1}{N},1 \right\} \; (i=1,...,\ell). 
\ee
The renormalized numbers $\rho_i:= N\lambda_i$ have a natural meaning of {\it occupation numbers} of the sites $a_i\in X$: they indicate for each site $i\in\{1,...,\ell\}$ by how many ``particles'' it is occupied when the system is in the state $\delta_{i_1...i_N}=\delta_{a_{i_1}}\otimes\cdots\otimes\delta_{a_{i_N}}$. In terms of the $\lambda_i$'s, the right hand side of \eqref{M1'} can be re-written as 
\be \label{M1''}
   M_1 S \delta_{i_1...i_N} = \sum_{i=1}^\ell \lambda_i \delta_i.
\ee
Clearly, equations \eqref{lambdai}, \eqref{M1''} imply the inclusion $``\subseteq"$ in a). To infer the reverse inclusion, note that any probability measure $\lambda\in\calPoneN(X)$ can be decomposed into $N$ Dirac measures of size $\tfrac{1}{N}$, i.e. $\lambda=\tfrac{1}{N}\sum_{k=1}^N\delta_{i_k}$ for some $i_1,...,i_N\in\{1,...,\ell\}$. Rearranging the $i_k$ in nondecreasing order and using 
\eqref{M1'} shows that $\lambda$ is the image under $M_1$ of the $N$-point measure $S\delta_{i_1...i_N}$.  To infer b) we denote $\mu:=M_2S\delta_{i_1...i_N}$ and $\mu_{ij}:=\mu(\{ (a_i,a_j) \})$, so that $\mu = \sum_{i,j=1}^\ell \mu_{ij}\delta_{ij}$, and distinguish two cases. 
\\[2mm]
{\bf Case 1:} $i\neq j$. In this case, by \eqref{M2'}
\be \label{muij}
   \mu_{ij} = \frac{1}{N(N-1)} \sharp \Bigl\{ \mbox{ pairs }(i_m,i_n) 
   \mbox{ with }i_m=i, \, i_n=j\Bigr\} 
            \, = \, \frac{1}{N(N-1)} \rho_i \rho_j = \frac{N}{N-1}\lambda_i\lambda_j.
\ee
{\bf Case 2:} $i=j$. In this case, \eqref{M2'} gives
$$
   \mu_{ii} = \frac{1}{N(N-1)} \rho_i (\rho_i-1) 
   = \frac{N}{N-1} \lambda_i^2 - \frac{1}{N-1}\lambda_i.
$$
Altogether we obtain that $\mu$ is given by the expression \eqref{extr}, establishing b). 
\\[2mm]
{\bf Proof of \eqref{supset}.} The finite set of pair states \eqref{extr} lies on the continuous manifold $\varphi_N(\calP(X))$, which -- geometrically -- is the image of a simplex under a continuous nonlinear map. Our strategy is to unearth, and use, the differential geometry of this manifold. 

\begin{figure}[http!]
\begin{center}
          \includegraphics[width=0.6\textwidth]{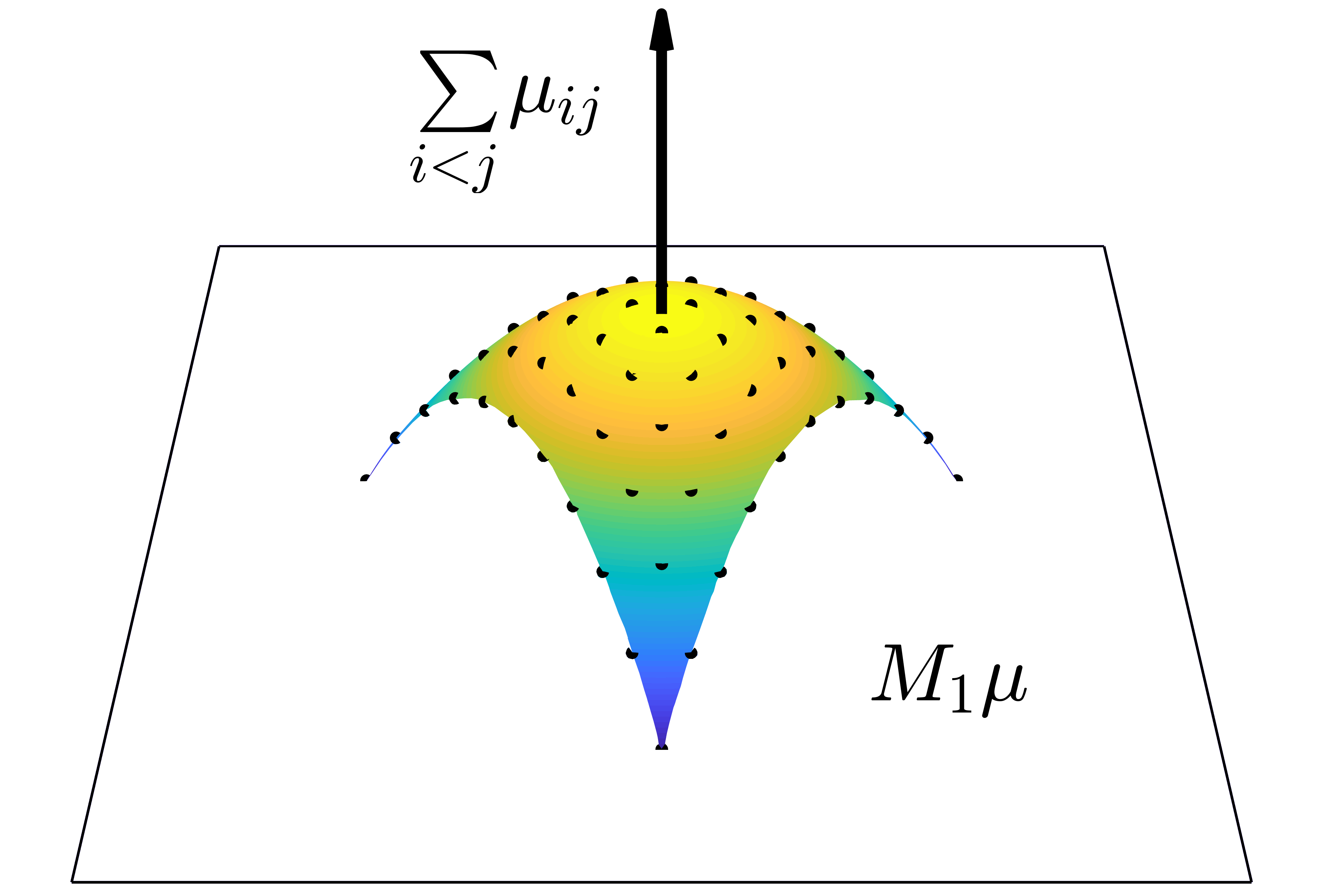}
\end{center}
\caption{The curved continuous manifold $\calM$ and the injective linear image of the set of extremal $N$-representable two-point measures. The picture shows the case $N=10$, $\ell=3$, in which there are 66 extreme points. The horizontal coordinates correspond to the one-point marginal. The vertical coordinate indicates the Wasserstein cost of the measure with respect to the discrete metric (see \eqref{wasser}, \eqref{metric}), which -- due to its linear dependence on the measure -- provides a coordinate direction in the space of two-point measures, and is useful for understanding the geometry of $N$-representable measures.} \label{F:manifold}
\end{figure}

More precisely, we restrict attention to a suitable lower-dimensional projection $\calM$ of this manifold. To construct this projection we introduce the following linear map $R \, : \, \calP_{sym}(X^2)\to \R^{\ell+1}$ which maps $\calP_{sym}(X^2)$ into a lower-dimensional vector space:  
\be \label{R}
   R\mu := \Bigl( M_1\mu, \frac{N-1}{N}\sum_{1\le i<j\le \ell}\mu_{ij}\Bigr) 
         = \Bigl( \sum_j\mu_{1j}, ..., \sum_j\mu_{\ell j}, 
            \frac{N-1}{N}\sum_{1\le i<j\le \ell}\mu_{ij}\Bigr). 
\ee 
We now introduce the continuous manifold $\calM := R \varphi_N(\calP(X))$ (see Figure \ref{F:manifold}), which by construction contains all points in $R \varphi_N(\calPoneN(X))$. By the explicit formula for $\varphi_N$ (see \eqref{UN}), we have $R \varphi_N(\lambda) = (\lambda, \sum_{1\le i < j \le \ell} \lambda_i\lambda_j)$, for any $\lambda\in\calP(X)$. It follows that $\calM$ has the structure of a graph of a scalar function:
\be \label{M}
   \calM = \{ (\lambda, g(\lambda)) \, : \, \lambda \in \calP(X) \}, 
   \mbox{ where }g(\lambda):=\sum_{1\le i<j\le \ell} \lambda_i\lambda_j.    
\ee
\noindent
The key to the extremality assertion we need to prove lies in the following differential-geometric fact: 
\begin{lemma}[Concavity] \label{L:conc} The function $g$ is strictly concave on $\calP(X)$; that is to say, 
\be \label{conc}
   g(\lambda') \; < \; g(\lambda) + \nabla g(\lambda)\cdot (\lambda'-\lambda) \mbox{ whenever }\lambda, \, \lambda'\in\calP(X), \, \lambda\neq\lambda'. 
\ee
\end{lemma}
{\bf Proof.} Since $g$ is quadratic, by Taylor expansion we have
\be \label{tay}
   g(\lambda') = g(\lambda) + \Bigl\langle \nabla g(\lambda), \, \lambda'-\lambda\Bigr\rangle 
   + \frac12 \Bigl\langle D^2g(\lambda)(\lambda'-\lambda), \, \lambda'-\lambda \Bigr\rangle.
\ee
Here $\langle a,b\rangle$ is the standard inner product $\sum_{i=1}^\ell a_i b_i$, and we employ the usual notation for gradient and Hessian, $(\nabla g)_i := \frac{\partial g}{\partial \lambda_i}$, $(D^2g)_{ij} = \frac{\partial^2g}{\partial \lambda_i\partial\lambda_j}$. For $g$ as defined in \eqref{M}, one calculates
\be \label{gradhess} 
   \nabla g(\lambda) = v - \begin{pmatrix} \lambda_1 \\ \vdots \\ \lambda_\ell \end{pmatrix},
   \;\;\;
   D^2g(\lambda) = v\, v^T - I, 
   \mbox{ with } v = \begin{pmatrix} 1 \\ \vdots \\ 1 \end{pmatrix}.  
\ee
Here $I$ denotes the identity matrix. 
It follows that $D^2g$ has eigenvalue $\ell-1$ on Span$\, v$, and eigenvalue $-1$ on the orthogonal complement $(\mbox{Span}\, v)^\perp = \{ \lambda\in\R^\ell \, : \, \sum_{i=1}^\ell \lambda_i=0\}$. This together with \eqref{tay} establishes \eqref{conc}, since the difference $\lambda'-\lambda$ belongs to $(\mbox{Span}\, v)^\perp$ (the tangent space of $\calP(X)$ at $\lambda$). 
\\[2mm]
Geometrically, the concavity of $g$ means that its graph lies strictly below the tangent space to the graph at $(\lambda,\, g(\lambda))$, as is evident from rewriting \eqref{conc} as
\be \label{conc'}
  \Bigl\langle \begin{pmatrix} \lambda' \\ g(\lambda') \end{pmatrix}, 
  \; n_\lambda \Bigr\rangle \; < \; \Bigl\langle \begin{pmatrix} \lambda \\ g(\lambda) \end{pmatrix},
  \; n_\lambda \Bigr\rangle \mbox{ for all }\begin{pmatrix} \lambda' \\ g(\lambda') \end{pmatrix} \in \calM \backslash \left\{ \begin{pmatrix} \lambda \\ g(\lambda) \end{pmatrix}\right\},
\ee  
where $n_\lambda$ is the ``upward'' normal to $\calM$ at $(\lambda,\, g(\lambda))$, explicitly:
$$
   n_\lambda = \frac{1}{\sqrt{1 + |\nabla g(\lambda)|^2}} \begin{pmatrix} -\nabla g(\lambda) \\ 1 \end{pmatrix}. 
$$
In particular, it follows from \eqref{conc'} that 
\be \label{conc''}
  \Bigl\langle R \varphi_N(\lambda'), \, n_\lambda\Bigr\rangle \; < \; 
  \Bigl\langle R \varphi_N(\lambda), \, n_\lambda \Bigr\rangle \mbox{ whenever }\lambda, \, \lambda'\in 
  \calPoneN(X), \, \lambda\neq \lambda'.
\ee
Thus each $R\varphi_N(\lambda)$, $\lambda\in\calPoneN(X)$, is an extreme point of $conv\bigl(R \varphi_N(\calPoneN(X)\bigr)$, and so -- by the linearity of $R$ -- each $\varphi_N(\lambda)$, $\lambda\in\calPoneN(X)$, is an extreme point of $conv\bigl(\varphi_N(\calPoneN(X))\bigr)$. This establishes \eqref{supset}, and completes the proof of Theorem \ref{T:extreme}. 
\\[2mm]
Together with \eqref{subset}, inequality \eqref{conc''} shows that any extremal $N$-representable two-point measure \eqref{extr} uniquely solves the problem of maximizing the linear functional $\langle R\mu,\, n_\lambda\rangle$ over $\mu\in\calP_{N-rep}(X^2)$. This problem was derived on purely differential-geometric grounds, yet it has an interesting optimal transport meaning. To see this we drop the inessential normalization factor from $n_\lambda$ and calculate using \eqref{gradhess}
\begin{eqnarray}
   \langle R\mu, \, \sqrt{1+|\nabla g(\lambda)|^2} n_\lambda \rangle & = &
   - \langle M_1\mu, \, \nabla g(\lambda)\rangle + \frac{N-1}{N} \sum_{1\le i<j\le \ell} \mu_{ij} \nonumber \\
   &=& - \sum_{i,j} \mu_{ij} + \sum_{i} \lambda_i \sum_j \mu_{ij} + \frac{N-1}{N} \sum_{i<j}\mu_{ij} \nonumber \\
   &=& \sum_{i,j=1}^\ell c_{ij}\mu_{ij}, \;\;\; c_{ij} := -1 + \frac{\lambda_i + \lambda_j}{2} + \frac{N-1}{2N}(1 - \delta_{ij}), \label{wassercalc}
\end{eqnarray}
where for the last equality sign we have used the symmetry of $\mu$. The crucial term $\sum_{i,j}(1-\delta_{ij})\mu_{ij}$ in the last expression is nothing but the {\it Wasserstein cost} 
\be \label{wasser}
   W[\mu] := \sum_{i,j=1}^\ell d(a_i,a_j)^p \mu_{ij} = \int_{X\times X} 
   d(x,y)^p d\mu(x,y) \;\;\; (1\le p < \infty)
\ee
with respect to the discrete metric on $X$, 
\be \label{metric} 
   d(x,y) := \begin{cases} 1 & x \neq y \\
                           0 & x = y. \end{cases}
\ee
Note that for the discrete metric, the Wasserstein cost is independent of $p$. We summarize this finding as follows. 
\begin{corollary}[Extremal two-point measures as maximizers of the Wasserstein cost] \label{C:wasser} $\;$ \\[1mm]
a) Each extreme point $\mu_*$ of $\calP_{N-rep}(X^2)$ is the unique solution of the problem
\be \label{J}
   \mbox{Maximize }J_\lambda[\mu] := \frac{2N}{N-1} \sum_i \lambda_i (M_1\mu)_i + W[\mu] 
   \mbox{ over }\calP_{N-rep}(X^2),
\ee
where $\lambda$ is the one-point marginal of $\mu_*$. \\[1mm]
b) In particular, each extreme point $\mu_*$ of $\calP_{N-rep}(X^2)$ is the unique solution of the optimal transport problem 
\be \label{J'}
   \mbox{Maximize }W[\mu] \mbox{ over }\mu\in\calP_{N-rep}(X^2) \mbox{ subject to }\mu\mapsto\lambda,
\ee
where $\lambda$ is the marginal of $\mu_*$. 
\end{corollary}
The variational problem \eqref{J} above -- which arose from purely differential-geometric considerations -- can be viewed as a ``soft-constraint'' version of the optimal transport problem \eqref{J'}, since the functional $J_\lambda$ promotes but does not rigidly enforce that the marginal $M_1\mu$ is close to $\lambda$. 
\section{Extremal symmetric $N$-point measures} 
\label{sec:extremal'} 
The extreme points of the set of symmetric $N$-point measures are straightforward to determine (see Proposition \ref{P:extreme'}). By comparing with the results of the previous section, we will see that this high-dimensional set has just the same number of extreme points as the much lower-dimensional set of $N$-representable two-point measures. We find this phenomenon quite remarkable. For instance it has the following optimal transport implication: every extreme point of the set of symmetric $N$-point measures, i.e. every unique optimizer of some OT problem \eqref{kant} with arbitrary $N$-body cost, must also be the unique optimizer of some OT problem with pairwise cost \eqref{pairwise}. In fact, by using the results of Section \ref{sec:extremal} we will see that a natural and universal such OT problem, maximizing the Gangbo-\'{S}wi\k{e}ch cost associated with the discrete metric on $X$ under a soft marginal constraint, does the job. See Corollary \ref{C:GS} below. 
\begin{proposition}[Extremal symmetric $N$-point measures] \label{P:extreme'} \textcolor{white}{.}
\\[1mm]
a) \emph{(Form of extreme points)} A probability measure on $X^N$ is an extreme point of $\calP_{sym}(X^N)$ if and only if it is a symmetrized Dirac measure, i.e. of form
\be \label{extr'}
    S \delta_{i_1...i_N} \mbox{ for some }i_1,...,i_N\in\{1,...,\ell\} \mbox{ with }i_1\le ... \le i_N, 
\ee
where $S$ is the symmetrization operator defined in \eqref{S}. 
\\[1mm]
b) \emph{(Parametrization by index vectors)} 
Different index vectors yield different extreme points; that is to say the map $(i_1,...,i_N)\mapsto S\delta_{i_1...i_N}$ maps the index set $\calI = \{(i_1,...,i_N)\in\{1,...,\ell\}^N\, : \, i_1\le ... \le i_N\}$ bijectively to the set of extreme points \eqref{extr'}. 
\end{proposition}
{\bf Proof.} First we show a). Let $E$ denote the set of probability measures of form \eqref{extr'}. 
We already found in Section \ref{sec:extremal} that $\calP_{sym}(X^N)$ is the convex hull of $E$ (see \eqref{MNS}). Hence the set of extreme points of $\calP_{sym}(X^N)$ is contained in $E$. So it suffices to show that any $\gamma\in E$ is not contained in the convex hull of $E\backslash\{\gamma\}$. This follows, for instance, because $\gamma=S\delta_{i_1...i_N}$ is the only element of $E$ whose support contains the point $(a_{i_1},...,a_{i_N})\in X^N$. 

It remains to show b). Injectivity of the map $(i_1,...,i_N)\mapsto S\delta_{i_1...i_N}$ from $\calI$ to $E$ follows from the above property of the support of $S\delta_{i_1,...,i_N}$. That the map is onto is clear from the definition of $E$. This proof of the proposition is complete.
\\[2mm]
The next result is less obvious. 
\begin{theorem}[Isomorphisms to sets of marginals] \label{T:extreme'}
\textcolor{white}{.} \\[1mm]
a) \emph{(Parametrization by one-point marginals)} The marginal map $M_1 \, : \, \calP_{sym}(X^N) \to \calP(X)$ maps the set of extreme points \eqref{extr'} bijectively to the whole range of $M_1$, that is to say to the set $\calPoneN(X)$ of $\frac{1}{N}$-quantized one-point measures. In particular, there are precisely 
$$
                    \tbinom{N+\ell-1}{\ell - 1}
$$
extreme points, and there exists an inverse map $\psi_N$ from $\calPoneN(X)$ to the set of extreme points such that for any $i_1,...,i_N$ as above, 
$$
     S\delta_{i_1...i_N} = \psi_N(M_1 S \delta_{i_1...i_N}). 
$$
b) \emph{(Relationship to extremal $N$-representable two-point measures)} The two-point marginal map $M_2 \, : \, \calP_{sym}(X^N)\to\calP_{N-rep}(X^2)$ is a bijection between the corresponding sets of extreme points. Moreover, for any extreme point $\mu$ of $\calP_{N-rep}(X^2)$, the extreme point of $\calP_{sym}(X^N)$ which is mapped to $\mu$ by $M_2$ is the unique representing measure of $\mu$ (see Def. \ref{D:Nrep}). 
\end{theorem}
Before coming to the proof of the theorem we note a corollary which relates the convex geometry of the set $\calP_{sym}(X^N)$ to the Gangbo-\'{S}wi\k{e}ch functional introduced in \cite{GS98} in the context of optimal transport on euclidean spaces, 
\be \label{GS}
         C_{GS}[\gamma] = \int_{X^N} \sum_{1\le i < j\le N} d(x_i,x_j)^p 
         d\gamma(x_1,...,x_N). 
\ee
In our setting of a finite state space $X$, we take $d$ to be the discrete metric \eqref{metric}. As for the ordinary Wasserstein cost \eqref{wasser} corresponding to $N=2$, the Gangbo-\'{S}wi\k{e}ch cost is trivially independent of $p$ when the underlying metric is the discrete metric.
\begin{corollary}[Extremal measures as maximizers of the Gangbo-\'{S}wi\k{e}ch cost] \label{C:GS} $\;$ Let $\gamma_*$ be any extreme point of $\calP_{sym}(X^N)$, and let $\lambda$ denote its one-point marginal. \\[1mm]
a) $\gamma_*$ is the unique solution of the problem 
\be \label{GS1}
   \mbox{Maximize }C_\lambda[\gamma] := N^2 \sum_i \lambda_i (M_1\gamma)_i 
   + C_{GS}[\gamma] \mbox{ over }\calP_{sym}(X^N).
\ee
b) In particular, $\gamma_*$ is the unique solution of the multi-marginal optimal transport problem  
\be \label{GS2}
   \mbox{Maximize the Gangbo-\'{S}wi\k{e}ch cost \eqref{GS} over }\gamma\in\calP_{sym}(X^N) 
   \mbox{ subject to }\gamma\mapsto\lambda.
\ee
\end{corollary}

Here, as before in Corollary \ref{C:wasser}, the variational problem in a) can be viewed as a soft-contraint version of the OT problem in b), in which closeness of the marginal $M_1\gamma$ to $\lambda$ is promoted but not enforced. 
\\[2mm]
{\bf Proof of Corollary \ref{C:GS}.} As for any pairwise cost, $C_{GS}[\gamma] = {N\choose 2} \int_{X\times X} d(x,y)^p d(M_2\gamma)(x,y)$, and so $\gamma$ is a solution to the variational problems in a) respectively b) if and only if it is a representing measure of a solution to the corresponding problems in Corollary \ref{C:wasser} b). The assertion now follows by combining the results of Corollary \ref{C:wasser} with the uniqueness result of Theorem \ref{T:extreme'} b). 
\\[2mm]
{\bf Proof of Theorem \ref{T:extreme'}.} We begin with a). As in the proof of Proposition 
\ref{P:extreme'} we denote the set \eqref{extr'} of symmetrized Dirac measures by $E$. 
That $M_1$ is onto from $E$ to $\calPoneN(X)$ was already shown in Lemma \ref{L:extreme} a). To show that $M_1$ is injective, we show that even the composed map from index vectors to one-point marginals, 
\be \label{composed}
   (i_1,...,i_N) \mapsto S\delta_{i_1...i_N} \mapsto M_1 S\delta_{i_1...i_N},
\ee
is injective on $\calI$. Indeed, if two index vectors $(i_1,...,i_N)$, $(i_1',...,i_N')\in\calI$ are different, then there must exist some $i\in\{1,...,\ell\}$ such that $\sharp\{k\in\{1,...,N\}\, : \, i_k=i\} \neq \sharp \{k\in\{1,...,N\} \, : \, i'_k = i\}$. By formulae \eqref{lambdaj}, \eqref{M1''} for the one-point marginal of a symmetrized Dirac measure, this implies that $M_1S\delta_{i_1...i_N}\neq M_1S\delta_{i_1'...i_N'}$. This establishes bijectivity and the existence of the map $\psi_N$. The assertion concerning the cardinality of $E$ now follows from the corresponding result for the cardinality of $\calPoneN(X)$ (see the proof of Corollary \ref{C:counting}). 
\\[2mm]
Now we deal with b). First we claim that $M_2E$ is equal to the set of extreme points of $\calP_{N-rep}(X^2)$. This follows by combining Lemma \ref{L:extreme}, which shows that $M_2 E$ is equal to the set of measures of form \eqref{extr}, and Theorem \ref{T:extreme}, which identifies the latter set as the set of extreme points of $\calP_{N-rep}(X^2)$. That $M_2$ is a bijection between $E$ and the set of extreme points of $\calP_{N-rep}(X^2)$ now follows because by Corollary \ref{C:counting} and a), both sets have the same cardinality. This esbalishes b), up to the uniqueness claim for representing measures. To establish the latter, assume that $\mu$ is any extreme point of $\calP_{N-rep}(X^2)$, and let $\gamma$ be any representing measure, i.e. any element of $\calP_{sym}(X^N)$ with $M_2\gamma = \mu$. By Proposition \ref{P:extreme'} a) and Minkowski's theorem, $\gamma$ belongs to the convex hull of $E$, that is to say 
\be \label{gammacomb} 
   \gamma = \sum_{i=1}^M \alpha_i\gamma_i \mbox{ for some }M\in\N, \; \gamma_i\in E, 
   \; \alpha_i\ge 0, \; \sum_i\alpha_i = 1. 
\ee
Applying the linear map $M_2$ to \eqref{gammacomb} yields $\mu = \sum_{i=1}^M \alpha_i M_2\gamma_i$. But by the injectivity of $M_2$ on $E$, the $M_2\gamma_i$ are all distinct. Since $\mu$ is by assumption an extreme point of $\calP_{N-rep}(X^2)$, it follows that $\alpha_{i_0}=1$ and $\alpha_i=0$ ($i\neq i_0$) for some $i_0$. Substitution of this finding into \eqref{gammacomb} shows $\gamma = \gamma_{i_0}$;
consequently $\gamma$ belongs to $E$. But uniqueness of representing measures $\gamma$ of $\mu$ within $E$ was already shown in the proof of the first part of b). The proof of Theorem \ref{T:extreme'} is complete.

\section{Measures with prescribed marginal} \label{sec:marginal}
Unconstrained linear optimization problems over the compact convex sets $\calP_{N-rep}(X^2)$ or $\calP_{sym}(X^N)$ will always attain their optimum at some extreme point, but multi-marginal optimal transport problems of the form \eqref{kant} or \eqref{redkant} typically will not. This is because most of the extreme points \eqref{extr} respectively \eqref{extr'} will fail to comply with a prescribed marginal condition. Probability measures in $\calP_{N-rep}(X^2)$ or $\calP_{sym}(X^N)$ with prescribed marginal, like any point inside a convex set, do however admit a (usually non-unique) representation as a convex combination of extreme points. To set the stage for the following section, here we write out this representation explicitly, and interpret it in the form of an integral over a suitable subset of $\calP(X)$, so as to bring out a formal analogy to the de Finetti representation for infinitely representable measures (see e.g. \cite{CFP15}). 

By Minkowski's theorem, for any $\gamma\in\calP_{sym}(X^N)$ respectively any $\mu\in\calP_{N-rep}(X^2)$ we have
\be \label{basic}
  \begin{array}{c}
        \mu = \sum\limits_{\lambda\in{\calPoneN(X)}} \alpha_\lambda
        \Bigl( \tfrac{N}{N-1} \lambda\! \otimes\! \lambda - 
    \tfrac{1}{N-1} \sum\limits_{i=1}^\ell \lambda_i \delta_i\otimes \delta_i  
    \Bigr) \mbox{ resp. }
        \gamma = \sum\limits_{\lambda\in\calPoneN(X)} \alpha_\lambda\psi_N(\lambda) \\
        \mbox{for some }\alpha_\lambda\ge 0 
        \mbox{ with }\sum\limits_{\lambda\in\calPoneN(X)}\alpha_\lambda = 1.
  \end{array}
\ee
Here $\psi_N$ is the bijective map from $\tfrac{1}{N}$-quantized probability measures on $X$ to symmetrized Dirac measures provided by Theorem \ref{T:extreme'}. 
The coefficients $\alpha_\lambda$ are highly non-unique in the case of $\mu$, but whenever $\gamma$ represents $\mu$ they can be taken to be the same in both expansions. This will be useful later. To reveal  an interesting analogy to de Finetti's theorem, we exploit the fact that the set of extreme points has been found to be isomorphic to a subset of $\calP(X)$, view the coefficients $(\alpha_\lambda)_{\lambda}$ as a probability measure $\alpha$ on this subset, i.e. a probability measure on $\calPoneN(X)$, and re-write \eqref{basic} as follows. Here and below the space of probability measures on $\calPoneN(X)$ will be denoted by $\calP(\calPoneN(X))$. 
\begin{lemma} \label{L:fancy} Given any symmetric $N$-point probability measures $\gamma$ on $X^N$,  respectively any $N$-representable two-point probability measure $\mu$ on $X^2$, there exists a probability measure $\alpha\in\calP(\calPoneN(X))$ such that
\be \label{fancy2}
   \gamma = \int_{\calPoneN(X)} \psi_N(\lambda) \, d\alpha(\lambda) 
\ee
respectively
\be \label{fancy1}
    \mu = \int_{\calPoneN(X)} \Bigl( \lambda\! \otimes\! \lambda + 
    \frac{1}{N-1}\bigl(\lambda\!\otimes\!\lambda - \sum_{i=1}^\ell \lambda_i \delta_i\otimes \delta_i  
    \bigr)\Bigr) \, d\alpha(\lambda),
\ee
where $\psi_N$ is the map provided by Theorem \ref{T:extreme'}.  
\end{lemma}
Note that in the limit $N\to\infty$ the integrand in \eqref{fancy1} tends to the independent measure $\lambda\otimes\lambda$ and the set $\calPoneN(X)$ of $\frac{1}{N}$-quantized probability measures tends -- formally -- to the set $\calP(X)$ of all probability measures, so \eqref{fancy1} turns into the celebrated de Finetti representation of infinitely representable measures on $X^2$,
\be \label{deF}
   \mu = \int_{\calP(X)} \lambda\!\otimes\!\lambda \, d\alpha(\lambda) 
   \mbox{ for some }\alpha \in \calP(\calP(X)).  
\ee
(For more information about \eqref{deF}, as well as an application to infinite-marginal optimal transport problems, we refer the reader to \cite{CFP15}.) 
In the representation \eqref{fancy2}--\eqref{fancy1}, any marginal constraint $\mu\mapsto\lambda_*$ or $\gamma\mapsto\lambda_*$ for some given $\lambda_*\in\calP(X)$ turns into the following constraint on the measure $\alpha$:
\be \label{margfancy}
   \int_{\calPoneN(X)} \lambda \, d\alpha(\lambda) = \lambda_*
\ee
or, in pedestrian notation, 
\be \label{margbasic}
   \sum_{\lambda\in\calPoneN(X)} \! \alpha_\lambda \, \lambda = \lambda_*.
\ee
The above representation formulae for $\mu$ and $\gamma$ as ``averages'' of extremal states have not, at this point, achieved any dimension reduction. The coefficient vector $(\alpha_\lambda)_{\lambda\in\calPoneN(X)}$ is still of length $\tbinom{N+\ell-1}{\ell - 1}$ --   exactly the same as the linear dimension of the original vector space of symmetric measures on $X^N$. In particular, in the typical situation in physical applications when the number of sites, $\ell$, is of the order of a constant times $N$ for some constant bigger than $1$ (i.e., the number of possible sites grows proportionally to the number of particles), the length of the coefficient vector increases exponentially as $N$ gets large.  
\section{Sparse averages of extremal states} \label{sec:sparse} 
We now show that the extreme points of the set of symmetric probability measures on $X^N$ with prescribed marginal $\lambda_*\in\calP(X)$, 
\be \label{Plambda*'}
    \calP_{sym,\lambda_*}(X^N) := \{ \gamma\in\calP_{sym}(X^N) \, : \, M_1\gamma = \lambda_*\}.
\ee
can be obtained by using only very sparse coefficient vectors in the expansion \eqref{basic}. The same will be established for the set of $N$-representable probability measures on $X^2$ with prescribed marginal,
\be \label{Plambda*}
    \calP_{N-rep,\lambda_*}(X^2) := \{ \mu\in\calP_{N-rep}(X^2) \, : \, M_1\mu = \lambda_*\}.
\ee 
We remark that the sets \eqref{Plambda*'} and \eqref{Plambda*} are always nonempty, a simple reason being that these spaces contain the $N$-fold respectively two-fold tensor product of the marginal with itself. 

In the following, for any convex set $K$ the set of its extreme points is denoted by $ext \, K$.
\begin{lemma}[Sparsity of extremal Kantorovich plans]\label{L:sparse} Let $\lambda_*\in\calP(X)$. 
\\[1mm]
a) Let $\calP_{SAE,\lambda_*}(X^N)$ denote the set of probability measures on $X^N$ of form \eqref{fancy2} for some $\alpha\in\calP(\calPoneN(X))$ which satisfies \eqref{margfancy} and is supported on at most $\ell$ elements of $\calPoneN(X)$. Then
\be \label{incl'}
    ext\bigl(\calP_{sym,\lambda_*}(X^N)\bigr) \subseteq \calP_{SAE,\lambda_*}(X^N) \subseteq 
    \calP_{sym,\lambda_*}(X^N). 
\ee
b) Let $\calP_{N,SAE,\lambda_*}(X^2)$ denote the set of probability measures on $X^2$ of form \eqref{fancy1} for some $\alpha\in\calP(\calPoneN(X))$ which satisfies \eqref{margfancy} and is supported on at most $\ell$ elements of $\calPoneN(X)$. Then 
\be \label{incl}
    ext\bigl(\calP_{N-rep,\lambda_*}(X^2)\bigr) \subseteq \calP_{N,SAE,\lambda_*}(X^2) \subseteq
    \calP_{N-rep,\lambda_*}(X^2).
\ee
\end{lemma}
The inclusions \eqref{incl} and \eqref{incl'} say that the sparse ansatz underlying the sets in the middle produces only probability measures in the sets on the right (i.e., admissible trial states in the Kantorovich problems \eqref{kant}, \eqref{redkant}), and contains all their extreme points. 

We propose to call measures in the middle set of \eqref{incl'} {\it sparse averages of extremal states}, or SAE states for short, and claim that they are precisely the same states introduced in a more elementary manner in the Introduction.

To see this we begin by eliminating the abstract map $\psi_N$. By definition, any element of $\calP_{SAE,\lambda_*}(X^N)$ is of the form $\sum_{\nu=1}^\ell \alpha^{(\nu)}\psi_N(\lambda^{(\nu)})$ for some $\tfrac{1}{N}$-quantized one-point probability measures $\lambda^{(\nu)}$ and some nonnegative $\alpha^{(\nu)}$ which sum to $1$. Each $\lambda^{(\nu)}$ can be decomposed further as  
\be \label{lambdanu}
    \lambda^{(\nu)} = \frac{1}{N} \sum_{k=1}^N \delta_{T_k(a_\nu )} 
\ee 
for some points $T_1(a_\nu),...,T_N(a_\nu)\in X$. 
The underlying physical picture is that any $\frac{1}{N}$-quantized probability measure arises by dropping $N$ points anywhere on the state space $X$ and encoding it by the location of the points, i.e. by the associated empirical measure. By definition of the map $\psi_N$ (see Section \ref{sec:extremal'}) and the construction of a symmetrized Dirac measure whose one-point marginal is a given element of $\calPoneN(X)$ (see the proof of Lemma \ref{L:extreme}) we have  
\be \label{psiNexplicit}
    \psi_N(\lambda^{(\nu)})=S\delta_{T_1(a_\nu)}\otimes\cdots\otimes\delta_{T_N(a_\nu)}
\ee
(note that because of the presence of the symmetrization operator, it is immaterial that the indices of the points $T_k(a_\nu)$ do not appear in nondecreasing order). 

Now observe something interesting.       
The possible values of the image points $T_k(a_\nu)$ in \eqref{lambdanu} are the sites in $X$, and hence there are $\ell$ possible values. But thanks to our sparsity lemma, $\nu$ runs only over just as many values. 
Hence the collection of image points $\{T_k(a_\nu)\, : \, \nu\in\{1,...,\ell\}, \, k\in\{1,...,N\}\}$ can be interpreted as $N$ maps $T_1,...,T_N\, : \, X\to X$; the measure $\lambda^{(\nu)}$ is then the value distribution of the ensemble $\{T_k\}_{k=1,...,N}$ of maps at the point $a_\nu$ (see Figure \ref{F:SAE} in the Introduction). Moreover the vector of weight coefficients $\alpha^{(\nu)}$ -- originally a probability measure on $\calPoneN(X)$ -- can be identified with a vector of ``site weights'', i.e. with a probability measure on the much smaller space $X$.   

To summarize:
\begin{proposition}[Characterization of SAE states]\label{P:SAE}
For any given marginal $\lambda_*\in\calP(X)$, and any probability measure on $X^N$, the following are equivalent: \\[1mm]
(i) The probability measure is of the form
\be \label{SAEelem}
         \gamma = \sum_{\nu=1}^\ell \alpha^{(\nu)} S \delta_{T_1(a_\nu)}\otimes\cdots\otimes
         \delta_{T_N(a_\nu)}
\ee 
for some maps $T_1,...,T_N\, : \, X\to X$ and some site weights
$\alpha^{(1)},...,\alpha^{(\ell)}\ge 0$ satisfying the system of equations
\be \label{constraint}
         \sum_{\nu=1}^\ell \lambda_i^{(\nu)}\alpha^{(\nu)} = (\lambda_*)_i \;\; (i=1,...,\ell),
\ee
where $\lambda^{(\nu)}$ is the empirical value distribution of the ensemble of maps $\{T_k\}_{k=1,...,N}$ at the site $a_\nu$, i.e. $\lambda^{(\nu)}$ is given by \eqref{lambdanu} or equivalently $\lambda^{(\nu)}_i=\tfrac{1}{N}\sharp\{k\in\{1,...,N\}\, : \, T_k(a_\nu)=a_i\}$. 
\\[1mm]
(ii) The probability measure is of the form 
\be \label{SAEadv}
     \gamma = \sum_{\nu=1}^\ell \alpha^{(\nu)} \psi_N(\lambda^{(\nu)})
\ee
for some probability measures $\lambda^{(1)},...,\lambda^{(\ell)}$ on $X$ which are $\tfrac{1}{N}$-quantized (i.e. belong to $\calPoneN(X)$) and some site weights $\alpha^{(1)},...,\alpha^{(\ell)}\ge 0$ satisfying the system of equations \eqref{constraint}, where $\psi_N$ is the isomorphism from $\calPoneN(X)$ to the set of symmetrized Dirac measures provided by Theorem \ref{T:extreme'}. 
\\[1mm]
Moreover if $\gamma$ is of form (i) then it is of form (ii) with the same $\lambda^{(\nu)}$s and $\alpha^{(\nu)}$s; in particular, it depends on the family of maps $\{T_k\}_{k=1,...,N}$ only through the associated value distributions \eqref{lambdanu}. Likewise, if $\gamma$ is of form (ii), then it is of the form (i) with the same $\lambda^{(\nu)}$s and $\alpha^{(\nu)}$s and with $T_1,...,T_N$ being {\it any} maps from $X$ to $X$ such that \eqref{lambdanu} holds.  
\end{proposition}
Note that the system of equations \eqref{constraint} automatically entails that the $\alpha^{(\nu)}$ sum to $1$, as is seen by summation over $i$.

As already pointed out in the Introduction, the number of parameters required to specify an SAE state is only $\ell\cdot N + \ell$, because each map $T_k$ is specified by the indices of its $\ell$ values at the sites $a_1,...,a_\ell$, and only $\ell$ coefficients $\alpha^{(\nu)}$ are required.

The meaning of the system \eqref{lambdaav} is that the {\it average} of the value distributions
$\lambda^{(\nu)}$ over the sites $a_\nu$ with respect to the 
site weights $\alpha^{(\nu)}$, 
\be \label{lambdaav}
    \lambda^{av} := \sum_{\nu=1}^\ell \alpha^{(\nu)} \lambda^{(\nu)},
\ee
must be equal to the prescribed marginal measure. By contrast, the standard multi-marginal Monge ansatz prescribes the distribution of values averaged over sites for the individual maps instead of the ensemble of maps.

For a graphical representation of a typical SAE state see the right panel in Figure \ref{F:SAE} in the Introduction. 

{\bf Proof of Lemma \ref{L:sparse}.} The second inclusion in \eqref{incl} and \eqref{incl'} is obvious from the representation formulae in Lemma \ref{L:fancy} for general 
$N$-representable two-point probability measures and general symmetric $N$-point probability measures. The key point is the first inclusion. We only prove it for \eqref{incl'}, the case of \eqref{incl} being analogous. Suppose $\gamma$ is an extreme point of $\calP_{sym,\lambda_*}(X^N)$. Like any other element of $\calP_{sym}(X^N)$, by Lemma \ref{L:fancy} it can be written in the form \eqref{fancy2} for some probability measure $\alpha\in\calP(\calPoneN(X))$. Suppose that the support of $\alpha$ contains more than $\ell$ elements, that is to say
$$
     \{\lambda \in \calPoneN(X)\, : \, \alpha_\lambda>0\} = \{ \lambda^{(1)},...,\lambda^{(\ell')} \}
     \mbox{ for some }\ell' \ge \ell+1.
$$
Clearly, using the shorthand notation $\alpha^{(\nu)}:=\alpha_{\lambda^{(\nu)}}$ we have
\be \label{expaninit}
    \gamma = \sum_{\nu=1}^{\ell'} \alpha^{(\nu)}\psi_N\bigl(\lambda^{(\nu)}\bigr). 
\ee
Moreover since $\gamma$ belongs to $\calP_{sym,\lambda_*}$, it satisfies the marginal condition \eqref{margbasic}, which now takes the following form
\be \label{inhom}
     \begin{pmatrix}
\lambda_1^{(1)}  &  \;\;\; \cdots \;\;\;   & \lambda_1^{(\ell')} \\
\vdots           &  \;\;\; \textcolor{white}{\cdots} \;\;  & \vdots \\
\lambda_{\ell}^{(1)} & \;\;\; \cdots \;\;\; & \lambda_{\ell}^{(\ell')}
     \end{pmatrix}
     \begin{pmatrix} 
\alpha^{(1)} \\[3mm]
\vdots \\[3mm]
\alpha^{(\ell')} 
    \end{pmatrix}
=  
    \begin{pmatrix}
(\lambda_*)_1 \\[1mm]
\vdots \\[1mm]
(\lambda_*)_{\ell}  
    \end{pmatrix}.  
\ee
This is a system of $\ell$ linear equations for $\ell'\ge \ell+1$ variables. Hence the corresponding homogeneous equation (\eqref{inhom} with the right hand side replaced by the zero vector) possesses a nonzero solution $\beta=(\beta^{(1)},...,\beta^{(\ell')})^T$. Moreover, by taking the sum of the $\ell$ homogeneous equations we see that 
$$
        0 = \sum_{i=1}^\ell \sum_{\nu=1}^{\ell'} \lambda_i^{(\nu)}\beta^{(\nu)} 
          = \sum_{\nu=1}^{\ell'} \beta^{(\nu)}.
$$
It follows that the vector $\beta$ contains at least one positive component and at least one negative component. Hence there exist $\eps_+,\, \eps_- >0$ such that 
\begin{eqnarray}
   & & \min_{\nu\in\{ 1,...,\ell'\} }\bigl(\alpha^{(\nu)} + \eps_+ \beta^{(\nu)}\bigr) = 0, 
       \label{min1} \\
   & & \min_{\nu\in\{ 1,...,\ell'\} }\bigl(\alpha^{(\nu)} - \eps_- \beta^{(\nu)}\bigr) = 0. 
       \label{min2}
\end{eqnarray} 
Define two new measures $\alpha_{\pm}\in\calP(\calPoneN(X))$ by 
\be \label{endpts}
   (\alpha_{\pm})_{\lambda} := \begin{cases} 
           \alpha^{(\nu)}\pm\eps_\pm \beta^{(\nu)} & \mbox{if }\lambda = \lambda^{(\nu)} \mbox{ for some }\nu\in\{1,...,\ell'\} \\
           0 & \mbox{otherwise.}
   \end{cases}          
\ee
It follows from the above definition that $\alpha$ is a convex combination of $\alpha^{\pm}$, more precisely 
\be \label{convalpha}
    \alpha = c_- \alpha_- + c_+\alpha_+ \; \mbox{ with }c_-=\tfrac{\eps_+}{\eps_+ + \eps_-}, \;\;             c_+ = \tfrac{\eps_-}{\eps_+ + \eps_-}.
\ee
Moreover by \eqref{min1} and \eqref{min2}, the support of $\alpha_{\pm}$ contains at most $\ell'-1$ points. Now introduce the following measures on $X^N$ 
$$
    \gamma_{\pm} := \int_{\calPoneN(X)} \psi_N(\lambda) \, d\alpha_{\pm}(\lambda)
    = \sum_{\nu=1}^{\ell'} \alpha^{(\nu)}_\pm \psi_N\bigl(\lambda^{(\nu)}\bigr), 
$$
with the shorthand notation $\alpha_\pm^{(\nu)}=(\alpha_\pm)_{\lambda^{(\nu)}}$. 
By construction the $\gamma_{\pm}$ still satisfy the marginal condition, i.e. belong to $\calP_{sym,\lambda_*}(X^N)$, and satisfy $\gamma=c_-\gamma_-+c_+\gamma_+$. 
Now there are two cases: either the $\gamma_{\pm}$ coincide with $\gamma$ or they do not. 
The second case cannot occur, because by assumption $\gamma$ is an extreme point of $\calP_{sym,\lambda_*}$. Thus $\gamma=\gamma_-=\gamma_+$, and so we have found a sparser representation of $\gamma$, by a measure in $\calP(\calPoneN(X))$ whose support contains at most $\ell'-1$ elements. Repeating the above construction leads to a representation of $\gamma$ by a measure whose support contains at most $\ell$ elements. 

Lemma \ref{L:sparse} immediately leads to the main result of this paper. 
\begin{theorem}[Breaking the curse of dimension in multi-marginal optimal transport]\label{T:OTnew}
\textcolor{white}{.} \\[1mm]
a) For any number $N\ge 2$ of marginals, any finite state space $X$, 
any cost function $c_N \, : \, X^N \to \R\cup\{+\infty\}$, and any prescribed marginal $\lambda_*\in\calP(X)$, the Kantorovich problem \eqref{kant} admits a solution which is an SAE state (see Definition \ref{D:SAE}). Moreover when $c_N$ is symmetric, this SAE state is also a minimizer of the Kantorovich cost $C[\gamma]$ over all of $\calP(X^N)$ subject to $\gamma\mapsto\lambda_*$. 
\\[1mm]
b) If, in addition, the cost function is of pairwise form, i.e. $c_N(x_1,...,x_N)=\sum_{1\le i<j\le N}c(x_i,x_j)$ for some $c\, : \, X^2\to\R\cup\{+\infty\}$, the minimum cost in \eqref{kant} is the same as that of the following explicit reduced problem: Minimize the functional
\be \label{red1}  
   I[\alpha,\lambda^{(1)},...,\lambda^{(\ell)}] = 
   \sum_{\nu=1}^\ell \alpha^{(\nu)}\Bigl(\tfrac{N^2}{2} \int_{X\times X} c(x,y)  d\lambda^{(\nu)}(x)
   d\lambda^{(\nu)}(y) - \tfrac{N}{2}\int_X c(x,x)\, d\lambda^{(\nu)}(x)
   \Bigr) 
\ee
over $\alpha=(\alpha^{(1)},...,\alpha^{(\ell)})\in\R^\ell$ and  $\lambda^{(1)},...,\lambda^{(\ell)}\in\calPoneN(X)$ subject to the constraints
\be
   \alpha^{(\nu)}\ge 0 \;\; (\nu=1,...,\ell), \;\;\;\;
   \sum_{\nu=1}^\ell\alpha^{(\nu)}\lambda^{(\nu)}=\lambda_*. 
\ee
Moreover if $(\alpha,\lambda^{(1)},...,\lambda^{(\ell)})$ is any minimizer of the reduced problem, then the associated SAE state \eqref{SAEadv} is a minimizer of \eqref{kant}. 
\end{theorem}
{\bf Proof of Theorem \ref{T:OTnew}.} By compactness and convexity of the set $\calP_{sym,\lambda_*}(X^N)$, the problem \eqref{kant} admits a solution which is an extreme point of this set. Assertion a) now follows from Lemma \ref{L:sparse} a). The proof of b) follows in an analogous manner, by using \eqref{elid} and Lemma \ref{L:sparse} b) together with the explicit representation \eqref{extr} for extremal $N$-representable measures.

\section{A new characterization of Monge states} \label{sec:Monge}
We now compare the sparse-average-of-extremals (SAE) ansatz with the classical but not always sufficient Monge ansatz
\be \label{Monge}
   \gamma = \sum_{\nu=1}^\ell \tfrac{1}{\ell} \delta_{T_1(a_\nu)}\otimes\cdots\otimes\delta_{T_N(a_\nu)} \mbox{ for $N$ permutations }T_1,...,T_N \, : \, X\to X. 
\ee  
(Recall that by re-ordering the sum one could without loss of generality assume $T_1=id$, but in the discussion below it will be convenient not to single out any particular map $T_k$.) The requirement that the $T_k$ be permutations implies that $\gamma$ has equal one-point marginals $\lambdabar$ (where $\lambdabar$ is the uniform measure \eqref{unif}). So the set of symmetrized Monge states,
\be \label{PMonge}
    \calP_{Monge}(X^N) = \{S\gamma \, : \, \gamma \mbox{ is of form \eqref{Monge}} \},
\ee
is contained in the set of symmetric $N$-point probability measures with uniform marginal. Hence it should be compared to the set of SAE states with uniform marginal, 
\be \label{PSAE}
    \calP_{SAE}(X^N) := \calP_{SAE,\lambdabar}(X^N) .
\ee 
\begin{theorem}[Characterization of Monge states] \label{T:Monge} A probability measure on $X^N$ is a symmetrized Monge state if and only if it is an SAE state with uniform marginal with all the site weights being equal to $1/\ell$, i.e. $\alpha^{(1)}=...=\alpha^{(\ell)}=1/\ell$. 
\end{theorem}
{\bf Proof.} The ``only if'' part is trivial, but the ``if'' part is not, because of the more stringent requirement in \eqref{Monge} that the $T_k$ must be permutations. To prove the ``if'' part, we start by taking any SAE state $\gamma$ with site weights $1/\ell$, that is to say
\be \label{start}
   \gamma = \sum_{\nu=1}^\ell \frac{1}{\ell} 
   \psi_N\bigl(\lambda^{(\nu)}\bigr)
\ee
for some $\lambda^{(1)},...,\lambda^{(\ell)}\in\calPoneN(X)$ which satisfy the system of equations
\be \label{constraint2}
   \sum_{\nu=1}^\ell \lambda_i^{(\nu)} \, \tfrac{1}{\ell} = \tfrac{1}{\ell} \;\; 
   (i=1,...,\ell). 
\ee
It is convenient here to use the abstract form \eqref{SAEadv} of these states, as the fact that $\gamma$ depends only on the $\lambda^{(\nu)}$ and $\alpha^{(\nu)}$ will be 
important. Consider the matrix whose columns are given by the $\lambda^{(\nu)}$s, i.e.
\be \label{bist}
   A_1 = \begin{pmatrix}
   \lambda_1^{(1)}  &  \cdots  & \lambda_1^{(\ell)} \\
   \vdots           &   & \vdots \\
   \lambda_{\ell}^{(1)} & \cdots  & \lambda_{\ell}^{(\ell)}
   \end{pmatrix}. 
\ee
By equations \eqref{constraint2} and the fact that the $\lambda^{(\nu)}$s are normalized, this matrix is doubly stochasic (i.e., it is a nonnegative matrix whose rows and columns sum to 1).  
Next we introduce a bipartite graph $G_1=(V_1 \cup V_1', E_1)$, as follows: the vertices are given by two disjoint sets of $\ell$ elements, say $V_1=\{1,...,\ell\}$, $V_1'=\{1',...,\ell'\}$, and the edge set is
\be \label{edges}
     E_1 = \bigl\{ \{i,j'\}\, : \, i,j\in V_1, \, \lambda_j^{(i)} > 0 \}. 
\ee
Thus replacing all nonzero entries of $A_1$ by $1$ would yield the adjacency matrix of the graph. 
By the Birkhoff-von Neumann theorem, $A_1$ can be written as a convex combination of permutation matrices. Therefore for any subset $W$ of $V_1$, the neighborhood $N(W)=\{v'\in V_1' \, : \, \{w,v'\}\in E_1 \mbox{ for some }w\in W\}$ satisfies
$$
    \sharp N(W) \ge \sharp W. 
$$
We now appeal to Hall's theorem for bipartite graphs (see e.g. \cite{GY05}) which implies that $G_1$ has a perfect matching, that is to say there exists a subset $M_1$ of $E_1$ such that each $v\in V_1$ belongs to some edge $\{v,w'\}\in M_1$ and different edges in $M_1$ are disjoint. Let $\tau_1 \, : \, V\to V$ be the map t
hat maps the element $v\in V$ to the element 
$(\tau_1(v))'$ in the same edge. Since $M_1$ is a perfect matching, $\tau_1$ is a permutation. Now subtract a suitable multiple of the corresponding permutation matrix from $A_1$, i.e. perform the following update:
\be \label{A2}
   A_2 := A_1 - \frac{1}{N} \sum_{i=1}^N e_{\tau_1(i)}e_i^T.
\ee
Since the components $(A_1)_{ij}$ of $A_1$ take values in $\{0,\tfrac{1}{N},...,\tfrac{N-1}{N},1\}$, it follows that $A_2\ge 0$. Now $\tfrac{N}{N-1}A_2$ is again doubly stochastic and by the same argumentation as before we can construct a permutation $\tau_2\, : \, \{1,...,\ell\}\to\{1,...,\ell\}$ such that $A_2 - \tfrac{1}{N}e_{T_2(i)}e_i^T \ge 0$. Iterating this procedure yields $N$ permutations $\tau_1,...,\tau_N$ such that 
\be\label{bist'}
  A_1 = \frac{1}{N}\sum_{k=1}^N \sum_{i=1}^N e_{\tau_k(i)}e_i^T. 
\ee
Comparing \eqref{bist} and \eqref{bist'} and passing from the permutations $\tau_k$ to the associated maps $T_k \, : \, X\to X$ defined by $T_k(a_\nu)=a_{\tau_k(\nu)}$ yields that the columns $\lambda^{(i)}$ of $A_1$ have the form $\lambda^{(i)}=\tfrac{1}{N}\sum_{k=1}^N\delta_{T_k(a_i)}$. The last part of Proposition \ref{P:SAE} now implies that $\gamma$ is equal to $\sum_{\nu=1}^\ell \tfrac{1}{\ell} S\, \delta_{T_1(a_\nu)}\otimes\cdots\otimes\delta_{T_N(a_\nu)}$ for the constructed maps $T_1,...,T_N$, i.e. it is a symmetrized Monge state.  
\\[5mm]
{\bf Acknowledgements.} We thank Maximilian Fichtl and S\"oren Behr for helpful discussions. 

\begin{small}

\end{small}

\end{document}